\documentclass[a4paper,12pt,final]{amsart}
\usepackage{times,a4wide,mathrsfs,amssymb,amsmath,amsthm,enumerate,xypic,tikzsymbols,dsfont}

\newcommand{\C}{\mathbb{C}}
\newcommand{\ZZ}{\mathbb{Z}}

\newcommand{\QQ}{\mathbb{Q}}
\newcommand{\NN}{\mathbb{N}}
\newcommand{\PP}{\mathbb{P}}
\newcommand{\LL}{\mathcal{L}}

\newcommand{\OO}{\mathcal O}

\newcommand{\XX}{\mathcal X}
\newcommand{\YY}{\mathcal Y}

\newcommand{\Zz}{\mathcal Z}
\newcommand{\EE}{\mathcal E}
\newcommand{\MM}{\mathcal M}

\newcommand{\pic}{\hbox{Pic}}

\newcommand{\wt}{\widetilde}
\newcommand{\rom}{\romannumeral}

\newcommand{\one}{\mathds{1}}

\DeclareMathOperator{\Pic}{Pic}
\DeclareMathOperator{\ima}{Im}

\DeclareMathOperator{\sym}{Sym}
\DeclareMathOperator{\Gr}{Gr}

\newtheorem{theorem}{Theorem}[section]
\newtheorem{claim}[theorem]{Claim}
\newtheorem{lemma}[theorem]{Lemma}

\newtheorem{corollary}[theorem]{Corollary}
\newtheorem{proposition}[theorem]{Proposition}

\newtheorem{remark}[theorem]{Remark}
\newtheorem{definition}[theorem]{Definition}
\newtheorem{convention}{Conventions}
\newtheorem{question}[theorem]{Question}
\newtheorem{notation}[theorem]{Notation}

\newtheorem{nonumbering}{Theorem}

\newtheorem{nonumberingc}{Corollary}

\newtheorem{nonumberingt}{Acknowledgements}

\begin{document}

\author[Robert Laterveer]
{Robert Laterveer}

\address{Institut de Recherche Math\'ematique Avanc\'ee,
CNRS -- Universit\'e 
de Strasbourg,\
7 Rue Ren\'e Des\-car\-tes, 67084 Strasbourg CEDEX,
FRANCE.}
\email{robert.laterveer@math.unistra.fr}

\title{Algebraic cycles and Fano threefolds of genus 8}

\begin{abstract} We show that prime Fano threefolds $Y$ of genus 8 have a multiplicative Chow--K\"unneth decomposition, in the sense of Shen--Vial. As a 
consequence, a certain tautological subring of the Chow ring of powers of $Y$ injects into cohomology.
 \end{abstract}

\thanks{Supported by ANR grant ANR-20-CE40-0023.}


\maketitle

\section{Introduction}

Given a smooth projective variety $Y$ over $\C$, let $A^i(Y):=CH^i(Y)_{\QQ}$ denote the Chow groups of $Y$ (i.e. the groups of codimension $i$ algebraic cycles on $Y$ with $\QQ$-coefficients, modulo rational equivalence). The intersection product defines a ring structure on $A^\ast(Y)=\bigoplus_i A^i(Y)$, the Chow ring of $Y$.
Motivated by the cases of K3 surfaces and abelian varieties, Beauville \cite{Beau3} has conjectured that for certain special varieties, the Chow ring should admit a multiplicative splitting. To make concrete sense of Beauville's elusive ``splitting property conjecture'', Shen--Vial \cite{SV} have introduced the concept of {\em multiplicative Chow--K\"unneth decomposition\/}. It seems both interesting and difficult to better understand the class of special varieties admitting such a decomposition (for some partial answers, cf. subsection \ref{ss:mck} below).

This paper aims to contribute to this program. The main result is as follows:

\begin{nonumbering}[=Theorem \ref{main}] Let $Y$ be a Fano threefold that is prime (i.e., the Picard group is generated by the canonical divisor) of genus 8. Then $Y$ has a multiplicative Chow--K\"unneth decomposition.
\end{nonumbering}

The argument proving Theorem \ref{main} relies on the classical fact that $Y$ is birational to a cubic threefold (cf. subsection \ref{ss:birat} below), plus the recently established fact that cubic threefolds have a multiplicative Chow--K\"unneth decomposition \cite{Diaz}, \cite{FLV2}.

We check that the multiplicative Chow--K\"unneth decomposition we construct behaves well with respect to the birationality with the cubic threefold (Proposition \ref{YandX}), and also with respect to birationalities among two genus 8 Fano threefolds (Proposition \ref{YandY}).

Using Theorem \ref{main}, we obtain a result concerning the {\em tautological ring\/}, which is a certain subring of the Chow ring of powers of $Y$:

\begin{nonumberingc}[=Corollary \ref{cor1}] Let $Y$ be a prime Fano threefold of genus 8, and $m\in\NN$. Let
  \[ R^\ast(Y^m):=\bigl\langle (p_i)^\ast(h), (p_{ij})^\ast(\Delta_Y)\bigr\rangle\ \subset\ \ \ A^\ast(Y^m)   \]
  be the $\QQ$-subalgebra generated by pullbacks of the polarization $h\in A^1(Y)$ and pullbacks of the diagonal $\Delta_Y\in A^3(Y\times Y)$. 
  The cycle class map induces injections
   \[ R^\ast(Y^m)\ \hookrightarrow\ H^\ast(Y^m,\QQ)\ \ \ \hbox{for\ all\ }m\in\NN\ .\]
   \end{nonumberingc}

 Corollary \ref{cor1} can be loosely paraphrased by saying that prime genus 8 Fano threefolds behave like hyperelliptic curves from the point of view of intersection theory (cf. Remark \ref{tava} below).  
   
To close this introduction, here is an open question:

\begin{question}\label{ques} Let $Y$ be a Fano threefold with $\Pic(Y)=\ZZ$. Does $Y$ have a multiplicative Chow--K\"unneth decomposition ? (And is Corollary \ref{cor1} true for $Y$ ?)
\end{question}

(Without the hypothesis on $\Pic(Y)$, the answer can be negative, as suggested by \cite[Example 1.7]{Beau3}.)   
 In the final section \ref{s: q}, we record another question on a related note.
   
 \vskip0.6cm

\begin{convention} In this article, the word {\sl variety\/} will refer to a reduced irreducible scheme of finite type over $\C$. A {\sl subvariety\/} is a (possibly reducible) reduced subscheme which is equidimensional. 

{\bf All Chow groups will be with rational coefficients}: we will denote by $A_j(Y)$ the Chow group of $j$-dimensional cycles on $Y$ with $\QQ$-coefficients; for $Y$ smooth of dimension $n$ the notations $A_j(Y)$ and $A^{n-j}(Y)$ are used interchangeably. 
The notation $A^j_{hom}(Y)$ will be used to indicate the subgroup of homologically trivial cycles.
For a morphism $f\colon X\to Y$, we will write $\Gamma_f\in A_\ast(X\times Y)$ for the graph of $f$.

The contravariant category of Chow motives (i.e., pure motives with respect to rational equivalence as in \cite{Sc}, \cite{MNP}) will be denoted 
$\MM_{\rm rat}$.
\end{convention}

\section{Preliminaries}

\subsection{MCK decomposition}
\label{ss:mck}

\begin{definition}[Murre \cite{Mur}] Let $X$ be a smooth projective variety of dimension $n$. We say that $X$ has a {\em CK decomposition\/} if there exists a decomposition of the diagonal
   \[ \Delta_X= \pi^0_X+ \pi^1_X+\cdots +\pi_X^{2n}\ \ \ \hbox{in}\ A^n(X\times X)\ ,\]
  such that the $\pi^i_X$ are mutually orthogonal idempotents and $(\pi_X^i)_\ast H^\ast(X,\QQ)= H^i(X,\QQ)$.
  
  (NB: ``CK decomposition'' is shorthand for ``Chow--K\"unneth decomposition''.)
\end{definition}

\begin{remark} The existence of a CK decomposition for any smooth projective variety is part of Murre's conjectures \cite{Mur}, \cite{J4}. 
\end{remark}

\begin{definition}[Shen--Vial \cite{SV}] Let $X$ be a smooth projective variety of dimension $n$. Let $\Delta_X^{sm}\in A^{2n}(X\times X\times X)$ be the class of the small diagonal
  \[ \Delta_X^{sm}:=\bigl\{ (x,x,x)\ \vert\ x\in X\bigr\}\ \subset\ X\times X\times X\ .\]
  An {\em MCK decomposition\/} is a CK decomposition $\{\pi_X^i\}$ of $X$ that is {\em multiplicative\/}, i.e. it satisfies
  \[ \pi_X^k\circ \Delta_X^{sm}\circ (\pi_X^i\times \pi_X^j)=0\ \ \ \hbox{in}\ A^{2n}(X\times X\times X)\ \ \ \hbox{for\ all\ }i+j\not=k\ .\]
  
 (NB: ``MCK decomposition'' is shorthand for ``multiplicative Chow--K\"unneth decomposition''.) 
  
  \end{definition}
  
  \begin{remark} The small diagonal (seen as a correspondence from $X\times X$ to $X$) induces the {\em multiplication morphism\/}
    \[ \Delta_X^{sm}\colon\ \  h(X)\otimes h(X)\ \to\ h(X)\ \ \ \hbox{in}\ \MM_{\rm rat}\ .\]
 Suppose $X$ has a CK decomposition
  \[ h(X)=\bigoplus_{i=0}^{2n} h^i(X)\ \ \ \hbox{in}\ \MM_{\rm rat}\ .\]
  By definition, this decomposition is multiplicative if for any $i,j$ the composition
  \[ h^i(X)\otimes h^j(X)\ \to\ h(X)\otimes h(X)\ \xrightarrow{\Delta_X^{sm}}\ h(X)\ \ \ \hbox{in}\ \MM_{\rm rat}\]
  factors through $h^{i+j}(X)$.
  
  If $X$ has an MCK decomposition, then setting
    \[ A^i_{(j)}(X):= (\pi_X^{2i-j})_\ast A^i(X) \ ,\]
    one obtains a bigraded ring structure on the Chow ring: that is, the intersection product sends $A^i_{(j)}(X)\otimes A^{i^\prime}_{(j^\prime)}(X) $ to  $A^{i+i^\prime}_{(j+j^\prime)}(X)$.
    
      It is expected that for any $X$ with an MCK decomposition, one has
    \[ A^i_{(j)}(X)\stackrel{??}{=}0\ \ \ \hbox{for}\ j<0\ ,\ \ \ A^i_{(0)}(X)\cap A^i_{hom}(X)\stackrel{??}{=}0\ ;\]
    this is related to Murre's conjectures B and D, that have been formulated for any CK decomposition \cite{Mur}.

  The property of having an MCK decomposition is severely restrictive, and is closely related to Beauville's ``splitting property' conjecture'' \cite{Beau3}. 
  To give an idea: hyperelliptic curves have an MCK decomposition \cite[Example 8.16]{SV}, but the very general curve of genus $\ge 3$ does not have an MCK decomposition \cite[Example 2.3]{FLV2}. For more detailed discussion, and examples of varieties with an MCK decomposition, we refer to \cite[Section 8]{SV}, as well as \cite{V6}, \cite{SV2}, \cite{FTV}, \cite{37},
  \cite{38}, \cite{39}, \cite{40}, \cite{46}, \cite{48}, \cite{55}, \cite{FLV2}, \cite{NOY}.
    \end{remark}

\subsection{Prime Fano threefolds of genus 8} The classification of Fano threefolds is one of the glories of 20th century algebraic geometry \cite{IP}. Fano threefolds that are {\em prime\/} (i.e. with Picard group of rank $1$ generated by the canonical divisor) come in 10 explicitly described families. In this paper we will be concerned with one of these families:

\begin{theorem}[Gushel \cite{Gu}, Mukai \cite{Mu}]\label{muk} Let $Y$ be a smooth projective Fano threefold with $\Pic(Y)=\ZZ$, of index $1$ and genus $8$. Then $Y$ is a dimensionally transverse intersection of the Grassmannian $\Gr(2,6)\subset\PP^{14}$ with a codimension $5$ linear subspace $\PP^9\subset\PP^{14}$.

Conversely, any smooth dimensionally transverse intersection $\Gr(2,6)\cap\PP^9\subset\PP^{14}$ is a prime Fano threefold of genus $8$.
\end{theorem}

\begin{proof} The last statement is just because $\Gr(2,6)$ is a Fano variety of dimension $8$ with Picard number 1, index $6$ and degree $14$; the codimension $5$ intersection $Y:=\Gr(2,6)\cap\PP^9$ thus has index $1$ and degree $d=14$ (i.e. genus $g=d/2+1=8$).
The first statement is proven in \cite{Gu} and independently in \cite{Mu}.
\end{proof}

\subsection{A birationality}\label{ss:birat} Any prime Fano threefold of genus $8$ is birational to a cubic threefold; this was already known to Fano \cite{Fa}, and has been used in many works \cite{IM1}, \cite{Isk}, \cite{Puts}, \cite{Ku}, \cite{Ku2}, \cite{KPS}.  The birationality is quite explicit:

\begin{theorem}[Puts \cite{Puts}, Iskovskikh \cite{Isk}]\label{puts} Let $Y$ be a prime Fano threefold of genus $8$. There exists a smooth cubic threefold $X$
and a birational map
  \[  \phi\colon\ \ Y\ \dashrightarrow\ X\ .\]
  More precisely, there exist a threefold $Z$ and morphisms
      \[  \xymatrix{ & Z \ar[dl]_q \ar[dr]^{p} \\
	Y && X ,}\]
	with the following properties:
	
\noindent
(\rom1) $p$ is the composition of the blow-up of a quintic elliptic curve $E_X\subset X$, and the blow-ups of (strict transforms of) lines $L_X^1,\ldots, L_X^{25}\subset X$; the lines $L_X^j$ are pairwise disjoint and $E_X\cap L_X^j$ consists of two points for each $j$;

\noindent
(\rom1) $q$ is the composition of the blow-up of a quintic elliptic curve $E_Y\subset Y$, and the blow-ups of (strict transforms of) lines $L_Y^1,\ldots, L_Y^{25}\subset Y$; the lines $L_Y^j$ are pairwise disjoint and $E_Y\cap L_Y^j$ consists of one point for each $j$;

\noindent
(\rom3) $p^{-1}(E_X)=q^{-1}(E_Y)$ and $p^{-1}(L_X^j)=q^{-1}(L_Y^j)$, $j=1,\ldots,25$.
 \end{theorem}
 
 \begin{proof} In a nutshell, the argument goes as follows (this can be found in \cite{Puts}, cf. also \cite[Chapter III Section 1]{Isk} and \cite[section 3]{IM1}). Let $V$ be a $6$-dimensional vector space, and let us write $G:=\Gr(2,V)\subset\PP(\wedge^2 V)$
 for the Grassmannian and its Pl\"ucker embedding. The projective dual to $G$ is the Pfaffian hypersurface $X_{13}\subset\PP(\wedge^2 V^\ast)$ (which is a cubic hypersurface of dimension $13$, singular along $\Gr(2,V^\ast)\subset X_{13}$). Given a smooth genus $8$ Fano threefold $Y=G\cap h_1\cap\cdots\cap h_5$, the intersection
   \[  X:= X_{13}\cap \langle h_1^\vee,\ldots,h_5^\vee\rangle\ \ \ \subset\ \PP(\wedge^2 V^\ast) \]
   is a smooth cubic threefold, called the {\em orthogonal cubic threefold\/}. A point $x\in X$ represents a skew-symmetric form of rank $4$ on $V$, whose kernel determines a line $n_x\subset\PP^5$.
   
   The {\em Palatini quartic\/} associated to $Y$ is defined as 
     \[ W_4:= \bigcup_{y\in Y} \ell_y\ \ \subset\ \PP^5\ ,\]
     where $\ell_y\subset\PP(V)\cong\PP^5$ is the line determined by $y\in Y$. This $W_4$ is a quartic hypersurface (singular along a curve), and it turns out that this coincides with
     \[ W_4= \bigcup_{x\in X} n_x\ \ \subset\ \PP^5\ ,\]
     where $X$ is the orthogonal cubic threefold. Choosing a general hyperplane $H\subset\PP^5$ and setting $W:=W_4\cap H$, one can define a map $u_Y\colon Y\to W$ by sending
     $y\in Y$ to the intersection point $\ell_y\cap H$. Likewise, one can define a map $u_X\colon X\to W$ by sending $x\in X$ to the intersection point
     $n_x\cap H$. One can check (cf. \cite{Puts}) that the maps $u_X$ and $u_Y$ are birational, and so one obtains a birationality
       \[   Y \ \xrightarrow{u_Y}\ W\ \xleftarrow{u_X}\ X\ .\]
    It remains to analyze the elimination of indeterminacy of this birational map. The indeterminacy locus of $u_Y$ (resp. $u_X$) is a quintic elliptic curve $E_Y$ (resp. $E_X$). For $H$ sufficiently general, the singular locus of $W$ consists of 25 singular points $w_1,\ldots,w_{25}$ (and these points correspond to the $L^j_Y$ via $u_Y$, and to the $L^j_X$ via $u_X$). Let $Z\to W$ be the blow-up of $w_1,\ldots,w_{25}$. One can show \cite[Section 6]{Puts} that there are morphisms from $Z$  
    to the blow-up of $Y$ in $E_Y$, and to the blow-up of $X$ in $E_X$; this gives the statements (\rom1), (\rom2), (\rom3).
    \end{proof}

  \subsection{An isomorphism of motives}
 
 \begin{proposition}\label{bla} Let $Y$ be a prime Fano threefold of genus $8$, and let $X$ be the orthogonal cubic threefold.
%
There is an isomorphism 
   \[ \ h(Y)\ \xrightarrow{\cong}\ h(X)\ \ \ \hbox{in}\ \MM_{\rm rat}\ .\] 
   \end{proposition}
 
 \begin{proof} Let $h\in A^1(Y)$ denote a hyperplane section (with respect to the embedding $Y\subset\PP^{14}$ given by Theorem \ref{muk}). We define
    \[   \begin{split}    \pi^0_Y&:= {1\over 14}\, h^3\times Y\ ,\\
                              \pi^2_Y&:= {1\over 14}\, h^2\times h\ ,\\
                               \pi^4_Y&:= {1\over 14}\, h\times h^2\ ,\\
                                \pi^6_Y&:= {1\over 14}\, Y\times h^3\ ,\\      
                                \pi^3_Y&:= \Delta_Y-\sum_{j\not=3} \pi^j_Y\ \ \ \in\ A^3(Y\times Y)\ ,\\
                                \end{split}\]
   and $h^j(Y):= (Y,\pi^j_Y,0)\in \MM_{\rm rat}$. Using a hyperplane section $h\in A^1(X)$, we similarly define a decomposition $h(X)=\oplus_j h^j(X)$.
 
 Let $\iota\colon D\hookrightarrow Z$ be the inclusion morphism of the exceptional divisor $D:= q^{-1}(E_Y)=p^{-1}(E_X)$ in $Z$ (where the curves
 $E_Y, E_X$ are as in Theorem \ref{puts}). The blow-up exact sequence gives equalities
   \begin{equation}\label{H3iso} \begin{split} \ima \Bigl(H^3(X,\QQ)\xrightarrow{p^\ast} H^3(Z,\QQ)\Bigr) &= \ker\Bigl( H^3(Z,\QQ)\to H^3(D,\QQ)\Bigr)\\
    &=  \ima \Bigl(H^3(Y,\QQ)\xrightarrow{q^\ast} H^3(Z,\QQ)\Bigr) \ ,\\
    \end{split} \end{equation}
and so in particular there is an isomorphism
  \[  p_\ast q^\ast\colon\ \ H^3(Y,\QQ)\ \xrightarrow{\cong}\ H^3(X,\QQ)\ .\]
  Using Manin's identity principle, this implies there is an isomorphism of homological motives
  \[  \Gamma_p\circ {}^t \Gamma_q\colon\ \ h^3(Y)\ \xrightarrow{\cong}\ h^3(X)\ \ \ \hbox{in}\ \MM_{\rm hom}\ .\]
It is known (cf. \cite[Theorem 4]{43}) that Fano threefolds have finite-dimensional motives (in the sense of Kimura \cite{Kim}), and so $h^3(X)$ and $h^3(Y)$ are finite-dimensional. It follows that there is also an isomorphism of Chow motives
   \[  \Gamma_p\circ {}^t \Gamma_q\colon\ \ h^3(Y)\ \xrightarrow{\cong}\ h^3(X)\ \ \ \hbox{in}\ \MM_{\rm rat}\ .\]      

For $j=0,\ldots,3$, we have $h^{2j}(Y)\cong h^{2j}(X)\cong\one(-j)$. Combined with the above isomorphism $h^3(Y)\cong h^3(X)$, this proves the proposition.
 \end{proof}

 \subsection{Franchetta property}
 
 \begin{definition} Let $\XX\to B$ be a smooth projective morphism, where $\XX, B$ are smooth quasi-projective varieties. We say that $\XX\to B$ has the {\em Franchetta property in codimension $j$\/} if the following holds: for every $\Gamma\in A^j(\XX)$ such that the restriction $\Gamma\vert_{X_b}$ is homologically trivial for the very general $b\in B$, the restriction $\Gamma\vert_b$ is zero in $A^j(X_b)$ for all $b\in B$.
 
 We say that $\XX\to B$ has the {\em Franchetta property\/} if $\XX\to B$ has the Franchetta property in codimension $j$ for all $j$.
 \end{definition}
 
 This property is studied in \cite{PSY}, \cite{BL}, \cite{FLV}, \cite{FLV3}.
 
 \subsection{Franchetta property for $Y\times Y$}
 
 \begin{notation}\label{not} Let $G$ be the Grassmannian $G:=\Gr(2,6)$, and let $\OO_G(1)$ be the polarization corresponding to the Pl\"ucker embedding
 $G\subset\PP^{14}$. Let 
   \[ B\ \subset\ \bar{B}:=\PP H^0(G,\OO_G(1)^{\oplus 5})\] 
   denote the Zariski open parametrizing smooth dimensionally transverse complete intersections, and let 
   \[ \YY\ \to\ B \]
   denote the universal family of smooth $3$-dimensional complete intersections (in view of Theorem \ref{muk}, this is the universal family of prime Fano threefolds of genus $8$).
  \end{notation}
 
 \begin{proposition}\label{Fr} Let $\YY\to B$ be the universal family of prime Fano threefolds of genus $8$. 
 The families $\YY\to B$ and $\YY\times_B \YY\to B$ have the Franchetta property.
  \end{proposition}
 
 \begin{proof} 
We exploit the relation with cubic threefolds (Theorem \ref{puts} and Proposition \ref{bla}). For a given $Y_b$, let $X_b$ denote the orthogonal cubic threefold, let $Z_b$ denote the resolution of indeterminacy with morphisms
    \begin{equation}\label{Zb}  \xymatrix{ & Z_b \ar[dl]_p \ar[dr]^{q} \\
	Y_b && X_b .}\end{equation}
      
 We observe that the construction of Theorem \ref{puts} makes sense familywise. That is, let us write $\XX\to B$ for the universal family of smooth dimensionally transverse cubic threefolds $X_b$ obtained as intersections $X_b=X_{13}\cap \PP^4_b$ (here a point $b\in B$ determines a $4$-dimensional linear subspace $\PP^4_b\subset\PP(\wedge^2 V^\ast)$ which is the orthogonal to $\PP^9_b\subset\PP(\wedge^2 V)$). Up to shrinking $B$, we may assume the $X_b$ avoid the singular locus of $X_{13}$, and thus $\XX$ is a smooth quasi-projective variety (it is an open in a projective bundle over $X_{13}\setminus\hbox{Sing}(X_{13})$).

 Fixing the choice of a hyperplane $H\subset V$, one has a birational map $\YY\dashrightarrow \XX$, and doing the required blow-ups, one obtains a diagram of $B$-schemes
   \[  \xymatrix{ & \Zz \ar[dl]_p \ar[dr]^{q} \\
	\YY && \XX ,}    \]
such that the base change along $b\hookrightarrow B$ gives back the diagram \eqref{Zb}. In particular, the isomorphism $h(Y)\cong h(X)$ of Proposition \ref{bla}
is generically defined, and so the induced isomorphisms of Chow groups
  \[  A^\ast(Y_b)\ \xrightarrow{\cong}\ A^\ast(X_b)\ ,\ \ \ A^\ast(Y_b\times Y_b)\ \xrightarrow{\cong}\ A^\ast(X_b\times X_b)\ \]
  send generically defined cycles to generically defined cycles. Consequently, Proposition \ref{Fr} reduces to proving the Franchetta property for the families $\XX\to B$ and $\XX\times_B \XX$. This is the content of the next proposition:

\begin{proposition}\label{Fr2} 	Let $\XX\to B$ be the universal family of orthogonal cubic threefolds, as above. 

\noindent
(\rom1)
Let $\Gamma\in A^{\ast}(\XX)$ be such that the restriction 
 $\Gamma\vert_{X_b}$ is homologically trivial, for the very general fiber $X_b$. Then 
   \[\Gamma\vert_{X_b}=0\ \ \ \hbox{in}\ A^\ast(X_b)\ \ \ \hbox{for\ all\ }b\in B\ .\] 
   
\noindent
(\rom2) Let $\Gamma\in A_\ast(\XX\times_B \XX)$ \footnote{Here we write $A_{\ast}(\XX\times_B \XX)$ rather than $A^\ast(\XX\times_B \XX)$ because $\XX\times_B \XX$ may not be smooth. The Gysin homomorphism $A_{\ast}(\XX\times_B \XX)\to A^\ast(X_b\times X_b)$ is still well-defined, because $B$ is smooth.} such that the restriction 
 $\Gamma\vert_{X_b\times X_b}$ is homologically trivial, for the very general fiber $X_b$. Then 
   \[\Gamma\vert_{X_b}=0\ \ \ \hbox{in}\ A^\ast(X_b\times X_b)\ \ \ \hbox{for\ all\ }b\in B\ .\] 
    \end{proposition}
 
 \begin{proof} (NB: note that this is not the standard presentation of cubic threefolds in $\PP^4$, and so this does not follow directly from the Franchetta property for the universal cubic in $\PP^4$ (and its square) !)
 
\noindent
(\rom1) Let us assume $j=2$ (since $A^1(X_b)=A^3(X_b)=\QQ$, the statement is non-trivial only for $j=2$).

Let $\bar{\XX}\to\bar{B}$ denote the Zariski closure. This is the universal family of all (possibly singular and degenerate) intersections of the Pfaffian $X_{13}\subset\PP^{14}$ with a $\PP^4$, and so $\bar{B}\cong\PP^r$. Since every point in $X_{13}$ imposes one condition on hyperplanes, the projection $\pi\colon\bar{\XX}\to X_{13}$ is a $\PP^{r-10}$-bundle.  
 
 Now given $\Gamma\in A_{\dim \XX-2}({\XX})$, let $\bar{\Gamma}\in A_{\dim \XX-2}(\bar{\XX})$ be a cycle restricting to $\Gamma$. The projective bundle
		formula yields
		\[  \bar{\Gamma}=  \pi^\ast(a_{2}) + \pi^\ast(a_1)\cdot \xi\ \ \  \hbox{in}\
		A_{\dim \XX-2}(\bar{\XX}) ,\]
		where $\xi\in A^1(\bar\XX)$ is relatively ample with respect to $\pi$ and
		$a_j\in A_{13-j}(X_{13})$.
		Let $h\in A^1(\bar{B})$ be a hyperplane section and let $\nu\colon \bar{\XX}
		\to \bar{B}$ denote the projection. We have
		\[  \nu^\ast(h)= c\, \xi+\pi^\ast(z)\ \ \ \hbox{in}\ \pic(\bar{\XX})\ ,\] 
		for some $c \in\QQ$ and $z\in A_{12}(X_{13})$. It
		is readily checked that the constant $c$ is non-zero. (Indeed, assume for a moment $c$
		were zero. Then we would have $\nu^\ast(h^{r})=\pi^\ast(z^{r})$
		in $H^{2r}(\bar{\XX})$. But the right-hand side is zero, since
		$r:=\dim\bar{B}>\dim X_{13}=13$, while the left-hand side is non-zero.)	
		
		The constant $c$ being non-zero, we can write
		\[ \xi= \pi^\ast(z)+\nu^\ast(d)\ \ \  \hbox{in}\  A_{\dim\XX-1}(\bar{\XX}) ,\]
		where $z\in A_{12}(X_{13})$ and $d\in A^1(\bar{B})$ are non-zero elements.
		The restriction of $\nu^\ast(d)$ to a fiber $X_b$ is zero,
		and so we find that
		\[  {\Gamma}\vert_{X_b}=  a_2^\prime\vert_{X_b}\ \ \  \hbox{in}\
		A^2(X_b) ,\]
		for some $a_2^\prime\in A_{11}(X_{13})$. (Note that the restriction map $A_{11}(X_{13})\to A^2(X_b)$ is well-defined, since $X_b\to X_{13}$ is the base change of $\PP^4\to\PP^{14}$ which is lci.)
		The proposition now follows from Lemma \ref{lem} below.
		\end{proof}
		
		\begin{lemma}\label{lem} $A_{11}(X_{13})\cong\QQ$.
		\end{lemma}
		
		To prove the lemma, we use the resolution of singularities constructed in 
		 \cite[Section 3]{Ku06}. Let $\wt{X}_{13}$ be the variety parametrizing pairs $(\omega,K)$, where $K\subset V$ is a $2$-dimensional subspace and 
		 $\omega$ is a skew-form containing $K$ in the kernel. The variety $\wt{X}_{13}$  is smooth (it is a projective bundle over a Grassmannian). The forgetful morphism $f\colon  \wt{X}_{13}\to X_{13}$ induces a diagram
		 \[ \begin{array}[c]{ccc}
		    E& \stackrel{i}{\hookrightarrow}& \wt{X}_{13}\\
		    &&\\
		   \ \  \ \ \ \downarrow{\scriptstyle f\vert_E}&&\ \ \downarrow{\scriptstyle f}\\
		    &&\\
		    S:=\hbox{Sing}(X_{13})&\hookrightarrow& \ \ X_{13}\ ,\\
		    \end{array}\]		 
where $f$ is an isomorphism over the non-singular locus $X_{13}^\circ:=X_{13}\setminus \hbox{Sing}(X_{13})$, and $f\vert_E\colon E\to S\cong\Gr(2,6)$ is a $\Gr(2,4)$-fibration.

The variety $\wt{X}_{13}$, being a projective bundle over a Grassmannian, has trivial Chow groups and so in particular the cycle class map induces an isomorphism
     \begin{equation}\label{ob4}   A^2(\wt{X}_{13})\ \xrightarrow{\cong}\ H^4(\wt{X}_{13},\QQ) \ .\end{equation}
The exceptional divisor $E$ has $ A_0(E)\cong\QQ $ (this follows from the fibration result in \cite[Theorem 3.2]{V5}, in view of the fact that both $S\cong\Gr(2,6)$ and the fibers $\Gr(2,4)$ of $f\vert_E$ have trivial Chow group of $0$-cycles). The Bloch--Srinivas argument \cite{BS} (plus the fact that $h^{pq}(E)=0$ for $p\not=q$) then implies that the cycle class maps induce isomorphisms
    \begin{equation}\label{ob5}   A^2(E)\ \xrightarrow{\cong}\ H^4(E,\QQ)\ .\end{equation}
     
    The isomorphisms \eqref{ob4} and \eqref{ob5} fit into a commutative diagram with exact rows 
    \begin{equation}\label{thisdiag} \begin{array}[c]{cccccc}
          A_{11}(E) &\xrightarrow{i_\ast}&  A_{11}(\wt{X}_{13}) &\to&  A_{11}(X_{13}) &\to 0\\
          &&&&&\\
         \  \downarrow{\cong} &&\ \downarrow{\cong}&&\downarrow&\\
          &&&&&\\
          H_{22}(E,\QQ) &\xrightarrow{i_\ast}&  H_{22}(\wt{X}_{13},\QQ) &\to&  W_{-22} H_{22}(X_{13},\QQ) &\to 0\\    
          \end{array}\end{equation}
     (Here vertical arrows are cycle class maps, and $W_\ast$ refers to Deligne's weight filtration on homology \cite{PS}.)  
     It follows that the right vertical arrow is also an isomorphism, and so one is reduced to proving that $  W_{-22} H_{22}(X_{13},\QQ) $ is one-dimensional.
     
     The exact sequence
     \[   \begin{array}[c]{cccccc}
     0 &\to& H^4(X_{13},\QQ) &\to& H^4(\wt{X}_{13},\QQ)\oplus H^4(S,\QQ)\ \to\ H^4(E,\QQ)\ \to\ &H^5(X_{13},\QQ)\\     
     &&\parallel &&&\parallel\\
      &&\QQ&&&0\\
      \end{array} \]
      (where the vertical equalities follow from weak Lefschetz) yields that
      \[  H_{22}(\wt{X}_{13},\QQ)\cong   H^4(\wt{X}_{13},\QQ)\cong H^4(X_{13},\QQ)\oplus H^2(S,\QQ)\cong\QQ^2\ .\]
      Since $f\vert_E\colon E\to S\cong\Gr(2,6)$ is a $\Gr(2,4)$-fibration, we find that
      \[   H_{22}(E,\QQ)\cong H^2(E,\QQ)\cong \QQ^2\ .\]
      The arrows labelled $i_\ast$ in \eqref{thisdiag} cannot be identically zero (indeed, let $h\in A^1(\wt{X}_{13})$ be an ample divisor. Then $i_\ast i^\ast(h)=h\cdot E$ is non-zero in $H_{22}(\wt{X}_{13},\QQ)$), and so we find that $A_{11}(X_{13})$ is one-dimensional, proving the lemma.
 
 \medskip     
 \noindent
 (\rom2) 
%
We consider again the Zariski closure $B\subset \bar{B}\cong\PP^r$, and the induced family $\bar{\XX}\times_{\bar{B}}\bar{\XX}\to\bar{B}$. Since 2 distinct points in $X_{13}$ impose 2 independent conditions on each hyperplane, the morphism
  \[ (\pi,\pi)\colon\ \ \bar{\XX}\times_{\bar{B}}\bar{\XX}\ \to\ X_{13}\times X_{13} \]
  has the structure of a $\PP^{r-20}$-bundle over $(X_{13}\times X_{13} )\setminus \Delta_{X_{13}}$, and a $\PP^{r-10}$-bundle over $  \Delta_{X_{13}}$.
  
  Using the projective bundle formula, this implies that for any $b\in B$ one has equality
    \begin{equation}\label{square} \begin{split} \ima\Bigl(  A_\ast(\bar{\XX}\times_{\bar{B}}\bar{\XX}\to A_\ast(X_b\times X_b)\Bigr)      = \ima\Bigl(  A_\ast(X_{13}\times X_{13})\to  A_\ast(X_b\times X_b)\Bigr)&\\  
    + \Delta_\ast \ima\Bigl( A_\ast(X_{13})\to &A_\ast(X_b)\Bigr)\ ,\\
    \end{split}\end{equation}
 where $\Delta\colon X_b\to X_b\times X_b$ denotes the diagonal embedding (the equality \eqref{square} is proven exactly as \cite[Proposition 4.1]{FLV}).
 As we have just seen (in the proof of (\rom1) directly above), there is equality
  \[\ima\Bigl( A_\ast(X_{13})\to A_\ast(X_b)\Bigr) = \langle h \rangle\ .\] 
  Moreover, the fact that $X_{13}$ has a resolution of singularities $\wt{X}_{13}$ with trivial Chow groups implies that $X_{13}$ has the Chow--K\"unneth property:
   \[ A_\ast(X_{13}\times X_{13})=A_\ast(X_{13})\otimes A_\ast(X_{13})\ .\]
   (This is just because $A_\ast(\wt{X}_{13})\to A_\ast(X_{13})$ is surjective, and $\wt{X}_{13}$, being a smooth projective variety with trivial Chow groups, has the Chow--K\"unneth property.)
   In view of these two observations, equality \eqref{square} reduces to
   \[   \ima\Bigl(  A_\ast(\bar{\XX}\times_{\bar{B}}\bar{\XX})\to A_\ast(X_b\times X_b)\Bigr)      =  \langle h, \Delta_{X_b} \rangle\ .\]
   It now remains to prove that for any smooth cubic threefold $X_b$ the $\QQ$-subalgebra 
     \[  R^\ast(X_b\times X_b):=\langle h, \Delta_{X_b} \rangle \ \  \subset\  A^\ast(X_b\times X_b)\]
   injects into cohomology via the cycle class map. In codimension $3$ this is fine, because the diagonal $\Delta_{X_b}$ is linearly independent from decomposable cycles
   $A^\ast(X_b)\otimes A^\ast(X_b)$ in $H^6(X_b\times X_b,\QQ)$ (indeed, if the diagonal were decomposable, all of $H^\ast(X_b,\QQ)$ would be algebraic, which is absurd since $h^{2,1}(X_b)\not=0$). For $j\not=3$, the equality
     \[  \Delta_{X_b}\cdot (p_j)^\ast(h) =\sum_{i=1}^4  {1\over 3}\,  h^i\times h^{4-i}\ \ \ \hbox{in}\ A^4(X_b\times X_b) \ \ \ (j=1,2)\]
 (which follows from the excess intersection formula, cf. \cite[Equation (12)]{FLV3}) implies that
  \[       R^j(X_b\times X_b)\ \ \subset\ R^\ast(X_b)\otimes R^\ast(X_b) = \langle h\rangle \otimes   \langle h\rangle\  \ \ (j\not=3)\ ,\]
  and hence $R^j(X_b\times X_b)\to H^{2j}(X_b\times X_b,\QQ)$ is injective also for $j\not=3$.  The proposition is proven.   
     \end{proof}
     
   For later use, we also record the following, which is an easy consequence of the above.  
     
 \begin{proposition}\label{FrYX} Let $\YY\to B$ and $\XX\to B$ be as above. The family $\YY\times_B \XX\to B$ has the Franchetta property.
 \end{proposition}
 
 \begin{proof} The (generically defined) isomorphism of Proposition \ref{bla} gives a commutative diagram
   \[ \begin{array}[c]{ccc}
            A^\ast(\YY\times_B \XX) & \to& A_\ast(\XX\times_B \XX)\\
              &&\\
              \downarrow&&\downarrow\\
              &&\\
              A^\ast(Y\times X) & \xrightarrow{\cong}& \ A^\ast(X\times X)\ .\\
              \end{array}\]
 The statement thus reduces to the Franchetta property for $\XX\times_B \XX$, which is Proposition \ref{Fr2}.             
    \end{proof}

 \subsection{Franchetta property for $Z\times Z$}
 
 \begin{notation}\label{notz} Let $\YY\to B$ be as in Notation \ref{not} the universal family of genus 8 prime Fano threefolds. 
Fixing a hyperplane $H\subset V$, one gets a family $\Zz\to B$ by blowing up $\YY$ along certain closed $B$-subschemes  (cf. proof of Proposition \ref{Fr}). 
(This family $\Zz\to B$ has the property that fiberwise one obtains the resolution of indeterminacy of Proposition \ref{puts}.)
 \end{notation}

 \begin{proposition}\label{Frz} Let $\Zz\to B$ be as in Notation \ref{notz}. The family $\Zz\times_B \Zz\to B$ has the Franchetta property in codimension $\le 3$.
  \end{proposition} 
  
  \begin{proof} 
 The $B$-scheme $\Zz$ is obtained from the $B$-scheme $\YY$ by blowing up
  the $B$-subscheme $\EE\to B$ (which is a family of elliptic curves), and then blowing up the strict transforms of the families of lines $\LL^i\to B$ ($i=1,\ldots, 25$). Let us write $Z:=Z_B$ for a fiber of $\Zz\to B$ (and likewise $E:=E_b$, $L^i:=L^i_b$ for a fiber of $\EE\to B$ resp. $\LL\to B$).
  The blow-up formula gives an isomorphism of motives (where arrows are the natural ones)
    \[ h(Z)\cong h(Y)\oplus h(E)(-1) \oplus \bigoplus_{i=1}^{25} \one(-1) \oplus \bigoplus_{i=1}^{25} \one(-2)  \ \ \ \hbox{in}\ \MM_{\rm rat}\ .\]
In particular, for Chow groups of the product we obtain isomorphisms (where arrows are the natural ones)
  \begin{equation}\label{isom} \begin{split}  A^j(Z\times Z)\cong  A^j(Y\times Y)\oplus A^{j-1}(Y\times E)\oplus A^{j-1}(E\times Y)\oplus   A^{j-2}(E\times E)&\\
        \oplus \bigoplus A^\ast(Y)\oplus \bigoplus A^\ast(E)\oplus &\bigoplus \QQ\ .\\
        \end{split}\end{equation}
    The isomorphism \eqref{isom} fits into a commutative diagram
    \[ \begin{array}[c]{ccccc cc}
          A^j(\Zz\times_B \Zz) &\!\! \xrightarrow{\phi} &   \!\! A^j(\YY\times_B \YY) & \oplus& \!\! \!\! A^{j-1}(\YY\times_B \EE)  \oplus    A^{j-1}(\EE\times_B \YY) &\!\! \oplus & \!\!A^{j-2}(\EE\times_B \EE)  \oplus \cdots\\         
       &&&&& &\\
       \downarrow && \downarrow && \downarrow  &&\downarrow  \\
       &&&&& &\\
       A^j(Z\times Z) &\!\! \xrightarrow{\cong} & \!\!   A^j(Y\times Y) & \oplus& \!\! \!\! A^{j-1}(Y\times E)  \oplus    A^{j-1}(E\times Y) & \!\!\oplus & \!\! A^{j-2}(E\times E)  \oplus \cdots\\  
       \end{array}\]
       (Here vertical arrows are fiberwise restrictions, and the arrow $\phi$ is defined because the isomorphism \eqref{isom} is defined explicitly in terms of the various blow-ups linking $Y$ to $Z$, and these blow-ups can be performed relatively.)
       The families $\YY\to B$ and $\YY\times_B \YY\to B$ have the Franchetta property (Proposition \ref{Fr}). Thus, in proving the proposition, we are reduced to proving the following two claims:
       
  \begin{claim}\label{c1} The families $\EE\to B$ and $\EE\times_B  \EE\to B$ have the Franchetta property in codimension $\le 1$.
  \end{claim}
  
  \begin{claim}\label{c2} The family $\YY\times_B \EE\to B$ has the Franchetta property in codimension $\le 2$.
  \end{claim}
     
 It remains to prove these two claims. As for Claim \ref{c1}, obviously it suffices to prove the Franchetta property for $\EE\times_B  \EE\to B$. Recall that
 $E:=E_Y\subset Y$ is the indeterminacy locus of the map $u_Y\colon Y\dashrightarrow W$ described in the proof of Theorem \ref{puts}.
  We exploit the specific description of $E$ given in \cite[Section 6]{Puts}: letting $H\subset V\cong\C^6$ denote the fixed hyperplane defining $W:=W_4\cap H$, the elliptic curve $E$ is given as
   \[ E= \Gr(2,H)\cap H_1\cap\cdots\cap H_5\ \ \ \subset\ \  \Gr(2,V)\cap H_1\cap\cdots\cap H_5=:Y\ .\]
This implies that $\EE\to B$ coincides with the universal complete intersection of $\Gr(2,H)\subset \PP^{14}$ with a codimension $5$ linear subspace. As in the proof of Proposition \ref{Fr}, the morphism
  \[ \bar{\EE}\times_{\bar{B}} \bar{\EE} \ \to\ \Gr(2,H)\times \Gr(2,H) \]
  has the structure of a $\PP^{r-10}$-bundle over $( \Gr(2,H)\times \Gr(2,H))\setminus \Delta_{\Gr(2,H)}$, and a $\PP^{r-5}$-bundle over $ \Delta_{\Gr(2,H)}$. Again as in the proof of Proposition \ref{Fr}, this implies
  \[ \ima   \Bigl( A^\ast(\bar{\EE}\times_{\bar{B}} \bar{\EE})\to A^\ast(E\times E)\Bigr) = \langle   h,\Delta_E\rangle\ =:R^\ast(E\times E)\ .\]
  It is readily checked that $R^j(E\times E)$ injects into cohomology for $j\le 1$ (it suffices to remark that $\Delta_E\in H^2(E\times E,\QQ)$ is not decomposable).
  This proves Claim \ref{c1}.
  
  As for Claim \ref{c2}, let $g\colon\Gr(2,H)\hookrightarrow \Gr(2,V)$ denote the inclusion morphism, and let $\Gamma_g$ denote its graph. The morphism
   \[ \bar{\EE}\times_{\bar{B}} \bar{\YY} \ \to\ \Gr(2,H)\times \Gr(2,V) \]
  has the structure of a $\PP^{r-10}$-bundle over $( \Gr(2,H)\times \Gr(2,V))\setminus \Gamma_g$, and a $\PP^{r-5}$-bundle over $\Gamma_g$. This is what is called a {\em stratified projective bundle\/} in \cite{FLV}. An application of \cite[Proposition 5.2]{FLV} then gives the inclusion
  \[  \begin{split}  \ima   \Bigl( A^j(\bar{\EE}\times_{\bar{B}} \bar{\YY})\to A^j(E\times Y)\Bigr) \ \ \subset\   \ima\Bigl( A^j(\Gr(2,H)\times \Gr(2,V))\to A^j(E\times Y)\Bigr) &\\
  + \ima\Bigl(  A^{j-3}(\Gamma_{g_b})\to A^j(E\times &Y)\Bigr)\ ,\\
      \end{split}\]
 where $g_b\colon E\hookrightarrow Y$ is the inclusion morphism.
 For $j\le 2$, the last summand obviously vanishes. Observing that 
   \[  A^\ast(  \Gr(2,H)\times \Gr(2,V))=A^\ast(\Gr(2,H))\otimes A^\ast(\Gr(2,V)) \ ,\] 
  and that the Franchetta property holds for $\EE\to B$ and for $\YY\to B$, this proves Claim \ref{c2}.   
    \end{proof}

  \section{Main result}
  
  \begin{theorem}\label{main} Let $Y$ be a prime Fano threefold of genus 8. Then $Y$ has a multiplicative Chow--K\"unneth decomposition.
  The induced bigrading on the Chow ring is such that
    \[ A^2_{(0)}(Y) =\QQ[c_2(Y)]\ .\]
  \end{theorem}  
  
\begin{proof} Let us first construct a CK decomposition for $Y$. Letting $h\in A^1(Y)$ denote a hyperplane section (with respect to the embedding $Y\subset\PP^{14}$ given by Theorem \ref{muk}), we set
    \[   \begin{split}    \pi^0_Y&:= {1\over 14}\, h^3\times Y\ ,\\
                              \pi^2_Y&:= {1\over 14}\, h^2\times h\ ,\\
                               \pi^4_Y&:= {1\over 14}\, h\times h^2\ ,\\
                                \pi^6_Y&:= {1\over 14}\, Y\times h^3\ ,\\      
                                \pi^3_Y&:= \Delta_Y-\sum_{j\not=3} \pi^j_Y\ \ \ \in\ A^3(Y\times Y)\ .\\
                                \end{split}\]
 We observe that this CK decomposition is {\em generically defined\/} with respect to the family $\YY\to B$ (Notation \ref{not}), i.e. it is obtained by restriction from
 ``universal projectors'' $\pi^j_\YY\in A^3(\YY\times_B \YY)$. (This is just because $h$ and $\Delta_Y$ are generically defined.)   
 
 Writing $h^j(Y):=(Y,\pi^j_Y)\in\MM_{\rm rat}$, we have
   \begin{equation}\label{h2j} h^{2j}(Y) \cong \one(-j)\ \ \ \hbox{in}\ \MM_{\rm rat}\ \ \ (j=0,\ldots, 3)\ .\end{equation}

 What we need to prove is that this CK decomposition is MCK, i.e.
      \begin{equation}\label{this} \pi_Y^k\circ \Delta_Y^{sm}\circ (\pi_Y^i\times \pi_Y^j)=0\ \ \ \hbox{in}\ A^{6}(Y\times Y\times Y)\ \ \ \hbox{for\ all\ }i+j\not=k\ ,\end{equation}
      or equivalently that
       \[   h^i(Y)\otimes h^j(Y)\ \xrightarrow{\Delta_Y^{sm}}\  h(Y) \]
       coincides with 
       \[ h^i(Y)\otimes h^j(Y)\ \xrightarrow{\Delta_Y^{sm}}\ h(Y)\ \to\ h^{i+j}(Y)  \ \to\ h(Y)\ , \]   
       for all $i,j$.
              
   As a first step, let us assume that we have three integers $(i,j,k)$ and at most one of them is equal to $3$. The cycle in \eqref{this} is generically defined and homologically trivial. The isomorphisms \eqref{h2j} induce an injection
   \[ (\pi^i_Y\times\pi^j_Y\times\pi^k_Y)_\ast A^6(Y\times Y\times Y)\ \hookrightarrow\ A^\ast(Y)\ ,\]
   and sends generically defined cycles to generically defined cycles (this is because the isomorphisms \eqref{h2j} are generically defined). As a consequence, the required vanishing \eqref{this} follows from the Franchetta property for $Y$, which is Proposition \ref{Fr}(\rom1).      
   
   In the second step, let us assume that among the three integers $(i,j,k)$, exactly two are equal to $3$. In this case, using the isomorphisms \eqref{h2j} we find an injection
   \[ (\pi^i_Y\times\pi^j_Y\times\pi^k_Y)_\ast A^6(Y\times Y\times Y)\ \hookrightarrow\ A^\ast(Y\times Y)\ ,\] 
respecting the generically defined cycles. As such,  
the required vanishing \eqref{this} follows from the Franchetta property for $Y\times Y$, which is Proposition \ref{Fr}(\rom2).

(As an aside, we note that in this second step, for most combinations of $(i,j,k)$ the vanishing \eqref{this} can also be proven by more elementary means;
only for $(i,j,k)=(3,3,4)$ do we need the full force of Proposition \ref{Fr}(\rom2).)
   
%
   
In the third and final step we treat the case $i=j=k=3$. For this case, we exploit the birationality with the cubic threefold $X$, in the guise of the following lemma:
   
   \begin{lemma}\label{pqequal} Let $Y$ be a prime Fano threefold of genus $8$, let $X$ be its orthogonal cubic threefold, and let
      \[  \xymatrix{ & Z \ar[dl]_q \ar[dr]^{p} \\
	Y && X }\] 
be as in Theorem \ref{puts}. 
There is equality 
  \[   {}^t \Gamma_p   \circ \pi^3_X\circ \Gamma_p =  {}^t \Gamma_q   \circ \pi^3_Y\circ \Gamma_q \ \ \ \hbox{in}\ A^3(Z\times Z)\ \]
  (where $\pi^3_X, \pi^3_Y$ are as above).
  \end{lemma}
  
 \begin{proof} Both sides of the equality are generically defined (with respect to the family $\Zz\to B$ of Proposition \ref{Frz}). The equality holds true in cohomology, because we have seen in equality \eqref{H3iso} that
   \[  \ima \Bigl(H^3(X,\QQ)\xrightarrow{p^\ast} H^3(Z,\QQ)\Bigr) 
    =  \ima \Bigl(H^3(Y,\QQ)\xrightarrow{q^\ast} H^3(Z,\QQ)\Bigr) \ .
    \]
The Franchetta property in codimension $3$ for $Z\times Z$ (Proposition \ref{Frz}) then allows to conclude that the equality of the lemma holds true modulo rational equivalence.
    \end{proof} 
    
 Armed with Lemma \ref{pqequal}, one can do the final step of the proof. The generically defined CK decomposition $\{\pi^j_X\}$ (cf. the proof of Proposition \ref{bla}) is an MCK decomposition (cf. \cite{FLV2}, or independently \cite{Diaz}). For brevity, let us write $m:=   {}^t \Gamma_p   \circ \pi^3_X\circ \Gamma_p  \in A^3(Z\times Z)$ for the correspondence of Lemma \ref{pqequal}. We note that ${}^t \pi^3_X=\pi^3_X$ and so ${}^t m=m$. We have equalities
   \[ \begin{split}  m\circ \Delta_Z^{sm}\circ (m\times m) & =  (m\times m\times m)_\ast (\Delta_Z^{sm}) \\
                                                                             & = \Bigl( ( {}^t \Gamma_p   \circ \pi^3_X\circ \Gamma_p)^{\times 3}\Bigr){}_\ast   (\Delta_Z^{sm}) \\
                                                                             & = (p,p,p)^\ast (\pi^3_X\times \pi^3_X\times \pi^3_X)_\ast (p,p,p)_\ast   (\Delta_Z^{sm}) \\
                                                                             & =  (p,p,p)^\ast (\pi^3_X\times \pi^3_X\times \pi^3_X)_\ast  (\Delta_X^{sm}) \\
                                                                             & =  (p,p,p)^\ast \Bigl(    \pi^3_X\circ \Delta_X^{sm}\circ ( \pi^3_X\times \pi^3_X)  \Bigr) \\
                                                                             & = 0\ \ \ \ \hbox{in}\ A^6(Z\times Z\times Z)\ .\\
                                                                             \end{split}\]
                              (Here the first and fifth equality are Lieberman's lemma \cite[Lemma 3.3]{V3}; the second equality is by definition of $m$; the third follows from general properties of correspondences; the fourth is the fact that $ (p,p,p)_\ast   (\Delta_Z^{sm})=\Delta_X^{sm}$, and the last equality expresses the fact that $\{\pi^j_X\}$ is an MCK decomposition.)         
                              
 On the other hand, we have equalities     
    \[ \begin{split}  \pi^3_Y\circ \Delta_Y^{sm}\circ ( \pi^3_Y\times \pi^3_Y) & = (\pi^3_Y\times \pi^3_Y\times \pi^3_Y)_\ast (\Delta_Y^{sm}) \\
                                                                                               & =  (q,q,q)_\ast (q,q,q)^\ast \Bigl(  (\pi^3_Y\times \pi^3_Y\times \pi^3_Y)_\ast (\Delta_Y^{sm}) \Bigr) \\
                                                                                               & = (q,q,q)_\ast   (q,q,q)^\ast   (\pi^3_Y\times \pi^3_Y\times \pi^3_Y)_\ast (q,q,q)_\ast (\Delta_Z^{sm})\\
                                                                                               & = (q,q,q)_\ast (m,m,m)_\ast (\Delta_Z^{sm})\\
                                                                                               & = 0 \ \ \ \ \hbox{in}\ A^6(Y\times Y\times Y)\ .\\
                                                                                               \end{split}\]     
This proves the required vanishing \eqref{this} for $i=j=k=3$, and shows that the CK decomposition is MCK.

Finally, let us prove that $A^2_{(0)}(Y)=\QQ[c_2(Y)]$. We remark that $c_2(Y)\in A^2_{(0)}(Y)$ because
   \[   (\pi^i_Y)_\ast c_2(Y) = (\pi^i_\YY)_\ast c_2(T_{\YY/B})\vert_Y = 0\ \ \ \hbox{in}\  A^2(Y)\ \ \  \hbox{for\ all\ } i\not= 4\ ,\]
   as follows from the Franchetta property for $\YY\to B$ (Proposition \ref{Fr}). One readily checks that $c_2(Y)$ is non-zero (e.g. one can take a smooth anticanonical section $S\subset Y$; if $c_2(Y)$ were zero then by adjunction also $c_2(S)=0$, which is absurd since $S$ is a K3 surface). Since $A^2_{(0)}(Y)$
   injects into $H^4(Y,\QQ)\cong\QQ$, it follows that $A^2_{(0)}(Y)=\QQ[c_2(Y)]$.   
    \end{proof}

 \section{Compatibilities} 
  The MCK decomposition for $Y$ behaves well with respect to the birationality of Theorem \ref{puts}:

 \begin{proposition}\label{YandX} Let $Y$ be a prime Fano threefold of genus 8, let $X$ be the orthogonal cubic threefold and let $\phi\colon Y \dashrightarrow X$
 be the birational map of Theorem \ref{puts}. One has
   \begin{equation}\label{part0}   \bar{\Gamma}_\phi\ \ \in\ A^3_{(0)}(Y\times X)\ ,\end{equation}
   where $Y\times X$ is given the product MCK decomposition.
   
   In particular, one has
     \[  \phi_\ast c_2(Y)  \in \QQ[c_2(X)]\ \ \ \subset\ A^2(X)\ .\]
     \end{proposition}
     
   \begin{proof} The point is that the map $\phi$, and hence also the cycle $ \bar{\Gamma}_\phi$, is generically defined (with respect to the base $B$).
   Since the product MCK decomposition is also generically defined, the result \eqref{part0} follows from the Franchetta property for $Y\times X$ (Proposition \ref{FrYX}).
   
   For the second part, we have just seen that $A^2_{(0)}(Y)=\QQ[c_2(Y)]$, and one also has $A^2_{(0)}(X)=\QQ[c_2(X)]$ (this is just because $c_2(X)$ is a sum of intersections of divisors, by adjunction). The result \eqref{part0} implies that
     \[ \phi_\ast  A^2_{(0)}(Y)  \ \ \subset\ A^2_{(0)}(X)\ \]
     (cf. \cite[Proposition 1.5]{SV2}), and so $\phi_\ast c_2(Y) \in \QQ[ c_2(X)]$. 
           \end{proof}
  
   The MCK decomposition for $Y$ also behaves well with respect to birational maps between genus 8 Fano threefolds:
   
   \begin{proposition}\label{YandY} Let $Y_1, Y_2$ be two prime Fano threefolds of genus 8 such that there is a birational map $\psi\colon Y_1\dashrightarrow Y_2$. 
   There is an isomorphism of motives
     \begin{equation}\label{motisom}   \bar{\Gamma}_\psi\colon\ \ h^3(Y_1)\ \xrightarrow{\cong}\ h^3(Y_2)\ \ \ \hbox{in}\ \MM_{\rm rat}\ ,\end{equation}
   and one has
             \begin{equation}\label{grade}   \bar{\Gamma}_\psi\ \ \in\ A^3_{(0)}(Y_1\times Y_2)\  \end{equation}
   (where $Y_1\times Y_2$ is given the product MCK decomposition).
   
   In particular,
          \[  \psi_\ast c_2(Y_1)  \in \QQ[c_2(Y_2)]\ \ \ \subset\ A^2(Y_2)\ .\]
   \end{proposition}
  
  \begin{proof} Let $X_1, X_2$ denote the cubic threefold orthogonal to $Y_1$ resp. $Y_2$, and let
  $\phi_i\colon Y_i\dashrightarrow X_i$ ($i=1,2$) denote the birational maps of Theorem \ref{puts}. Thanks to the Torelli theorem for cubic threefolds \cite{CG}, \cite{Tyu} (cf. also \cite[Remark 2.19]{Ku}), $X_1$ and $X_2$ are isomorphic, and there is a commutative diagram
    \[ \begin{array}[c]{ccc}
         Y_1 & \stackrel{\phi_1}{\dashrightarrow} & X_1 \\
         &&\\
        \ \ \  \downarrow{\scriptstyle \psi} &&  \ \  \ \ \downarrow{\scriptstyle \psi_X}  \\
        &&\\
           Y_2 & \stackrel{\phi_2}{\dashrightarrow} & X_2 \\
           \end{array}\]
     where $\psi_X$ is an isomorphism. Since the maps $\phi_i$ induce isomorphisms $h^3(Y_i)\cong h^3(X_i)$ (proof of Proposition \ref{bla}), this proves \eqref{motisom}.
     
 Using Manin's identity principle (applied to the diagram obtained from the above diagram by taking the cartesian product with a fixed variety), it follows that there is equality of cycle classes
      \[    \bar{\Gamma}_\psi=  {}^t \bar{\Gamma}_{\phi_2}\circ \Gamma_{\psi_X}\circ  \bar{\Gamma}_{\phi_1}\  \ \ \hbox{in}\ A^3_{}(Y_1\times Y_2)\  .\]    
The $\bar{\Gamma}_{\phi_i}$ are in $A^3_{(0)}(Y_i\times X_i)$ (Proposition \ref{YandX}), and composition of correspondences preserves the subring $A^\ast_{(0)}()$ \cite[Proposition 1.5]{SV2}. To conclude the truth of equation \eqref{grade}, it only remains to establish the following:

\begin{lemma} $\Gamma_{\psi_X}\in A^3_{(0)}(X_1\times X_2)$.
\end{lemma}

To prove the lemma, we observe that the isomorphism $\psi_X$ preserves (the canonical divisor and hence) the hyperplane class, which implies that
  \[  \pi^j_{X_2}= (\psi_X\times\psi_X)_\ast \pi^j_{X_1}= \Gamma_{\psi_X}\circ \pi^j_{X_1}\circ {}^t \Gamma_{\psi_X}\ \ \ \hbox{in}\ A^3(X_2\times X_2)\ \ \ 
          \forall j\ .\]
Consequently, we find that
  \[  \pi^k_{X_2}\circ \Gamma_{\psi_X}\circ \pi^j_{X_1}=0\ \ \ \forall\ k\not=j\ .\]
  Now writing
    \[ \begin{split} \Gamma_{\psi_X}&= \Bigl(\sum_{k} \pi^k_{X_2}\Bigr)\circ \Gamma_{\psi_X}\circ  \Bigl(\sum_{j} \pi^j_{X_1}\Bigr)\\  
                                                            &= \sum_j   \pi^j_{X_2}\circ \Gamma_{\psi_X}\circ  \pi^j_{X_1}      \\
                                                            &= \sum_j   (\pi_{X_1}^{6-j}\times \pi_{X_2}^j)_\ast ( \Gamma_{\psi_X})\ \ \ \hbox{in}\ A^3(X_1\times X_2)\ ,\\
                                                   \end{split}\]   
                                                   we find that $ \Gamma_{\psi_X}\in A^3_{(0)}(X_1\times X_2)$ as requested. This closes the proof of the lemma and the proposition.  
     \end{proof}
  
  \begin{remark}\label{kuz} Let $Y_1, Y_2$ be 2 genus 8 Fano threefolds that are birational, as in Proposition \ref{YandY}.
  Using the relation with the orthogonal cubic threefold, Kuznetsov proves \cite[Corollary 3.18]{Ku} that there is an equivalence of categories
    \[\hbox{Kuz}(Y_1)\cong \hbox{Kuz}(Y_2)\] 
(where $\hbox{Kuz}(Y_i)$ denotes the Kuznetsov component of the derived category of $Y_i$). It follows from \cite[Theorem 1.1]{B+} that the converse is also true: an equivalence of categories $\hbox{Kuz}(Y_1)\cong \hbox{Kuz}(Y_2)$ implies that $Y_1$ and $Y_2$ are birational.  
  \end{remark}

 \section{The tautological ring}
 
 \begin{corollary}\label{cor1} Let $Y$ be a prime Fano threefold of genus 8, and let $m\in\NN$. Let
  \[ R^\ast(Y^m):=\bigl\langle (p_i)^\ast(h), (p_{ij})^\ast(\Delta_Y)\bigr\rangle\ \subset\ \ \ A^\ast(Y^m)   \]
  be the $\QQ$-subalgebra generated by pullbacks of the polarization $h\in A^1(Y)$ and pullbacks of the diagonal $\Delta_Y\in A^3(Y\times Y)$. (Here $p_i$ and $p_{ij}$ denote the various projections from $Y^m$ to $Y$ resp. to $Y\times Y$).
  The cycle class map induces injections
   \[ R^\ast(Y^m)\ \hookrightarrow\ H^\ast(Y^m,\QQ)\ \ \ \hbox{for\ all\ }m\in\NN\ .\]
   \end{corollary}

\begin{proof} This is inspired by the analogous result for cubic hypersurfaces \cite[Section 2.3]{FLV3}, which in turn is inspired by analogous results for hyperelliptic curves \cite{Ta2}, \cite{Ta} (cf. Remark \ref{tava} below) and for K3 surfaces \cite{Yin}.

As in \cite[Section 2.3]{FLV3}, let us write $o:={1\over 14} h^3\in A^3(Y)$, and
  \[ \tau:= \Delta_Y - {1\over 14}\, \sum_{j=0}^3  h^j\times h^{3-j}\ \ \in\ A^3(Y\times Y) \]
  (this cycle $\tau$ is nothing but the projector on the motive $h^3(Y)$ considered above).
Moreover, let us write 
  \[ \begin{split}   o_i&:= (p_i)^\ast(o)\ \ \in\ A^3(Y^m)\ ,\\
                        h_i&:=(p_i)^\ast(h)\ \ \in \ A^1(Y^m)\ ,\\
                         \tau_{ij}&:=(p_{ij})^\ast(\tau)\ \ \in\ A^3(Y^m)\ .\\
                         \end{split}\]
We define the $\QQ$-subalgebra
  \[ \bar{R}^\ast(Y^m):=\langle o_i, h_i, \tau_{ij}\rangle\ \ \ \subset\ H^\ast(Y^m,\QQ) \]
  (where $i$ ranges over $1\le i\le m$, and $1\le i<j\le m$). One can prove (just as \cite[Lemma 2.12]{FLV3} and \cite[Lemma 2.3]{Yin}) that the $\QQ$-algebra $ \bar{R}^\ast(Y^m)$
  is isomorphic to the free graded $\QQ$-algebra generated by $o_i,h_i,\tau_{ij}$, modulo the following relations:
    \begin{equation}\label{E:X'}
			o_i\cdot o_i = 0, \quad h_i \cdot o_i = 0,  \quad 
			h_i^3 =14\,o_i\,;
			\end{equation}
			\begin{equation}\label{E:X2'}
			\tau_{i,j} \cdot o_i = 0 ,\quad \tau_{i,j} \cdot h_i = 0, \quad \tau_{i,j} \cdot \tau_{i,j} = -10\, o_i\cdot o_j
			\,;
			\end{equation}
			\begin{equation}\label{E:X3'}
			\tau_{i,j} \cdot \tau_{i,k} = \tau_{j,k} \cdot o_i\,;
			\end{equation}
			\begin{equation}\label{E:X4'}
			\sum_{\sigma \in \mathfrak{S}_{12}}  \prod_{i=1}^{6} \tau_{\sigma(2i-1), \sigma(2i)} = 0\,. 
			\end{equation}
(NB: with reference to the notation of \cite[Lemma 2.12]{FLV3}, we have $b_{pr}=\dim H^3(Y,\QQ)=10$.)

To prove Corollary \ref{cor1}, it remains to check that all these relations are already verified modulo rational equivalence.
The relations \eqref{E:X'} take place in $R^\ast(Y)$ and so they follow from the Franchetta property for $Y$ (Proposition \ref{Fr}). 
The relations \eqref{E:X2'} take place in $R^\ast(Y^2)$. The first and the last relations are trivially verified, because ($Y$ being Fano) $A^6(Y^2)=\QQ$. As for the second relation of \eqref{E:X2'}, this follows from the Franchetta property for $Y\times Y$ (Proposition \ref{Fr}(\rom2)). (Alternatively, one can deduce the second relation from the MCK decomposition: the product $\tau_{} \cdot h_i$ lies in $A^4_{(0)}(Y^2)$, and it is readily checked that $A^4_{(0)}(Y^2)$ injects into $H^8(Y^2,\QQ)$.)
   
   Relation \eqref{E:X3'} takes place in $R^\ast(Y^3)$ and follows from the MCK relation. Indeed, we have
   \[  \Delta_Y^{sm}\circ (\pi^3_Y\times\pi^3_Y)=   \pi^6_Y\circ \Delta_Y^{sm}\circ (\pi^3_Y\times\pi^3_Y)  \ \ \ \hbox{in}\ A^6(Y^3)\ ,\]
   which (using Lieberman's lemma) translates into
   \[ (\pi^3_Y\times \pi^3_Y\times\Delta_Y)_\ast    \Delta_Y^{sm}  =   ( \pi^3_Y\times \pi^3_Y\times\pi^6_Y)_\ast \Delta_Y^{sm}                                                         \ \ \ \hbox{in}\ A^6(Y^3)\ ,\]
   which means that
   \[  \tau_{13}\cdot \tau_{23}= \tau_{12}\cdot o_3\ \ \ \hbox{in}\ A^6(Y^3)\ .\]
   
  Finally, relation \eqref{E:X4'}, which takes place in $R^\ast(Y^{12})$, 
  expresses the fact that 
     \[ \sym^{12} \dim H^3(Y,\QQ)=0\ ,\]
     where $H^3(Y,\QQ)$ is viewed as a super vector space.
 To check that this relation is also verified modulo rational equivalence (i.e. relation \eqref{E:X4'} holds true in $A^{36}(Y^{12})$), we observe that relation \eqref{E:X4'} involves a cycle  
 contained in
   \[ A^\ast(\sym^{12} (h^3(Y))\ .\]
 But the Chow motive ($\sym^{11} h^3(Y)$ and a fortiori) $\sym^{12} h^3(Y)$ is zero because $\dim H^3(Y,\QQ)=10$ and $h^3(Y)$ is oddly
  finite-dimensional (all Fano threefolds have finite-dimensional motive \cite[Theorem 4]{43}). This ends the proof.
%
%
%
 \end{proof}

\begin{remark}\label{tava} Given a curve $C$ and an integer $m\in\NN$, one can define the {\em tautological ring\/}
  \[ R^\ast(C^m):=  \langle  (p_i)^\ast(K_C),(p_{ij})^\ast(\Delta_C)\rangle\ \ \ \subset\ A^\ast(C^m) \]
  (where $p_i, p_{ij}$ denote the various projections from $C^m$ to $C$ resp. $C\times C$).
  Tavakol has proven \cite[Corollary 6.4]{Ta} that if $C$ is a hyperelliptic curve, the cycle class map induces injections
    \[  R^\ast(C^m)\ \hookrightarrow\ H^\ast(C^m,\QQ)\ \ \ \hbox{for\ all\ }m\in\NN\ .\]
   On the other hand, there are many (non hyperelliptic) curves for which the tautological ring $R^\ast(C^3)$ does {\em not\/} inject into cohomology (this is related to the non-vanishing of the Ceresa cycle, cf. \cite[Remark 4.2]{Ta} and also \cite[Example 2.3 and Remark 2.4]{FLV2}). 
   
Corollary \ref{cor1} shows that genus 8 Fano threefolds behave similarly to hyperelliptic curves. It would be interesting to understand what happens for other  Fano threefolds: is Corollary \ref{cor1} true for all of them or not ? (This is related to Question \ref{ques}.)
\end{remark}

\section{Question}
\label{s: q}

\begin{question} Let $Y$ be a prime Fano threefold of genus $8$, and let $S\subset Y$ be a smooth anticanonical divisor. Then $S$ is a K3 surface (and the general genus 8 K3 surface arises in this way). Is it true that
  \[ \ima\bigl(  A^1(S)\ \to\ A^2(Y)\bigr) = \QQ[h^2]\ \ ??\]
 (To prove this, it would suffice to prove that the graph of the inclusion morphism $S\hookrightarrow Y$ is in $A^3_{(0)}(S\times Y)$, with respect to the product MCK decomposition. I have not been able to do this.)
 
 This question also makes sense for cubic threefolds and their anticanonical divisors. 
  \end{question}

 \vskip1cm
\begin{nonumberingt} Thanks to Charles Vial and the referee for helpful comments. This note is dedicated to Y.: you know why...\Winkey
\end{nonumberingt}


\vskip1cm
\begin{thebibliography}{dlPG99}

		

\bibitem{Beau3} A. Beauville, On the splitting of the Bloch--Beilinson filtration, in: Algebraic cycles and motives (J. Nagel and C. Peters, editors), London Math. Soc. Lecture Notes 344, Cambridge University Press 2007,

\bibitem{BV} A. Beauville and C. Voisin, On the Chow ring of a K3 surface, J. Alg. Geom. 13 (2004), 417---426,

\bibitem{BL} N. Bergeron and Z. Li, {Tautological classes on moduli space of
			hyperk\"ahler manifolds}, Duke Math. J., arXiv:1703.04733,
			
\bibitem{B+} M. Bernardara, E. Macr\`\i, S. Mehrotra and P. Stellari, A categorical invariant for cubic threefolds, Adv. Math. 229 no 2 (2012), 770---803,			

\bibitem{BS} S. Bloch and V. Srinivas, Remarks on correspondences and algebraic cycles, American Journal of Mathematics Vol. 105, No 5 (1983), 1235---1253,

\bibitem{CG} H. Clemens and P. Griffiths, The intermediate Jacobian of the cubic threefold, Annals of Math. 95 (1972),
281---356,

\bibitem{Diaz} H. Diaz, The Chow ring of a cubic hypersurface, International Math. Research Notices,


\bibitem{Fa} G. Fano, Sulle sezione spaziale della variet\`a grassmanniana delle rette dello spazio a cinque dimensioni, Rend.
Acc. Naz. Lincei 11 (1930), 329---335,

\bibitem{FLV} L. Fu, R. Laterveer and Ch. Vial, The generalized Franchetta conjecture for some hyper-K\"ahler varieties (with an appendix joint with M.
Shen), Journal Math. Pures et Appliqu\'ees (9) 130 (2019), 1---35,

\bibitem{FLV3} L. Fu, R. Laterveer and Ch. Vial, The generalized Franchetta conjecture for some hyper-K\"ahler varieties, II, Journal de l'Ecole Polytechnique--Math\'ematiques 8 (2021), 1065---1097,

\bibitem{FLV2} L. Fu, R. Laterveer and Ch. Vial, Multiplicative Chow--K\"unneth decompositions and varieties of cohomological K3 type, Annali Mat. Pura ed Applicata, doi: 10.1007/s10231-021-01070-0,

\bibitem{FTV} L. Fu, Z. Tian and Ch. Vial, Motivic hyperk\"ahler resolution conjecture for generalized Kummer varieties, arXiv:1608.04968,

\bibitem{F} W. Fulton, Intersection theory, Springer--Verlag Ergebnisse der Mathematik, Berlin Heidelberg New York Tokyo 1984,

\bibitem{Gu} N. Gushel, On Fano threefolds of genus 8, St. Petersburg Math. J. 4 no 1 (1993), 115---129,

\bibitem{IM1} A. Iliev and L. Manivel, Pfaffian lines and vector bundles on Fano threefolds of genus 8, J. Alg. Geom. 16 (2007), 499---530,

\bibitem{Isk} V. Iskovskikh, Birational automorphisms of three-dimensional algebraic varieties, J. Sov. Math. 13 (1980), 815---868,

\bibitem{IP} V. Iskovskih and Yu. Prokhorov, Algebraic Geometry V: Fano varieties, Encyclopaedia of Math. Sciences 47, Springer-Verlag, Berlin 1999,

\bibitem{J2} U. Jannsen, Motivic sheaves and filtrations on Chow groups, in: Motives (U. Jannsen et alii, eds.), Proceedings of Symposia in Pure Mathematics Vol. 55 (1994), Part 1,  

\bibitem{J4} U. Jannsen, On finite-dimensional motives and Murre's conjecture, in: Algebraic cycles and motives (J. Nagel and C. Peters, editors), Cambridge University Press, Cambridge 2007,

\bibitem{Kim} S.-I.~Kimura, Chow groups are finite dimensional, in some sense, Math. Ann. 331 no 1 (2005), 173---201, 
		
\bibitem{Ku} A. Kuznetsov, Derived categories of cubic and $V_{14}$ threefolds, Proc. V.A. Steklov Inst. Math 246 (2004), 183---207,

\bibitem{Ku06} A. Kuznetsov, Homological projective duality for Grassmannians of lines, arXiv:0610957,

\bibitem{Ku2} A. Kuznetsov, Derived categories of Fano threefolds, Proc. V. A. Steklov Inst. Math 264 (2009), 110---122,

\bibitem{KPS} A. Kuznetsov, Y. Prokhorov and C. Shramov, Hilbert schemes of lines and conics and automorphism groups of Fano threefolds, Japanese Journal of Mathematics 13 no 1 (2018), 109---185,



\bibitem{37} R. Laterveer, A remark on the Chow ring of K\"uchle fourfolds of type $d3$, Bulletin Australian Math. Soc. 100 no. 3 (2019), 410---418,

\bibitem{38} R. Laterveer and Ch. Vial, On the Chow ring of Cynk--Hulek Calabi--Yau varieties and Schreieder varieties, Canadian Journal of Math. 72 no 2 (2020), 505---536,

\bibitem{39} R. Laterveer, Algebraic cycles and Verra fourfolds, Tohoku Math. J. 72 no. 3 (2020), 451---485,

\bibitem{40} R. Laterveer, On the Chow ring of certain Fano fourfolds, Ann. Univ. Paedagog. Crac. Stud. Math. 19 (2020), 39---52,

\bibitem{46} R. Laterveer, Algebraic cycles and Gushel--Mukai fivefolds, Journal of Pure and Applied Algebra 225 no. 5 (2021), 
doi: 10.1016/j.jpaa.2020.106582,

\bibitem{48} R. Laterveer, Algebraic cycles and intersections of 2 quadrics, Mediterranean Journal of Mathematics (2021), 
doi: 10.1007/s00009-021-01787-5,

\bibitem{55} R. Laterveer, Algebraic cycles and intersections of a quadric and a cubic, Forum Mathematicum  33 no. 3 (2021), 845---855, 

\bibitem{Mu} S. Mukai, New developments in the theory of Fano threefolds: vector bundle method and moduli problems, Sugaku Expositions 15 no 2 (2002), 125---150,

\bibitem{Mur} J. Murre, On a conjectural filtration on the Chow groups of an algebraic variety, parts I and II, Indag. Math. 4 (1993), 177---201,

\bibitem{MNP} J. Murre, J. Nagel and C. Peters, Lectures on the theory of pure motives, Amer. Math. Soc. University Lecture Series 61, Providence 2013,

\bibitem{NOY} A. Negut, G. Oberdieck and Q. Yin, Motivic decompositions for the Hilbert scheme of points of a K3 surface, arXiv:1912.09320v1,

\bibitem{PSY} N. Pavic, J. Shen and Q. Yin, On O'Grady's generalized Franchetta conjecture, Int. Math. Res. Notices (2016), 1---13,

\bibitem{PS} C. Peters and J. Steenbrink, Mixed Hodge Theory, Ergebnisse der Math. 52, Springer Verlag 2008,

\bibitem{Puts} P. Puts, On some Fano-threefolds that are sections of Grassmannians, Indag. Math. 44 (1982), 77---90,

\bibitem{Sc} T. Scholl, Classical motives, in: Motives (U. Jannsen et alii, eds.), Proceedings of Symposia in Pure Mathematics Vol. 55 (1994), Part 1, 

\bibitem{SV} M. Shen and Ch. Vial, The Fourier transform for certain hyperK\"ahler fourfolds, Memoirs of the AMS 240 (2016), no.1139,

\bibitem{SV2} M. Shen and Ch. Vial, The motive of the Hilbert cube $X^{[3]}$, Forum Math. Sigma 4 (2016), 55 pp., 

\bibitem{Ta2} M. Tavakol, The tautological ring of the moduli space $M_2^{rt}$, International Math.
Research Notices 178 (2013), 

\bibitem{Ta} M. Tavakol, Tautological classes on the moduli space of hyperelliptic curves with rational tails, J. Pure Applied Algebra 222 no 8 (2018),
2040---2062,

\bibitem{Tyu} A. Tyurin, Geometry of the Fano surface of a nonsingular cubic $F\subset \PP^4$ and Torelli theorems for Fano
surfaces and cubics, Izv. USSR Math. 35 (1971), 498---529,

\bibitem{43} Ch. Vial, Projectors on the intermediate algebraic Jacobians, New York J. Math. 19 (2013), 793---822,

\bibitem{V5} Ch. Vial, Algebraic cycles and fibrations, Documenta Math. 18 (2013), 1521---1553,

\bibitem{V3} Ch. Vial, Remarks on motives of abelian type, Tohoku Math. J. 69 (2017), 195---220,

\bibitem{V6} Ch. Vial, On the motive of some hyperk\"ahler varieties, J. Reine Angew. Math. 725 (2017), 235---247,



\bibitem{Yin} Q.~Yin, Finite-dimensionality and cycles on powers of K3 surfaces, Comment. Math. Helv. 90 (2015), 503--511.


\end{thebibliography}
\end{document}